\documentclass[11pt,fleqn]{article}
\usepackage{amssymb,latexsym,amsmath,amsfonts}
\usepackage{graphicx}

    \advance\textheight by \topskip

    \numberwithin{equation}{section}

    \newtheorem{theorem}{Theorem}[section]
    \newtheorem{lemma}[theorem]{Lemma}

    \newtheorem{Definition}[theorem]{Definition}
    
    \newtheorem{Remark}[theorem]{Remark}
    \newenvironment{remark}{\begin{Remark}\rm}{\end{Remark}}
    \newtheorem{Example}[theorem]{Example}

    \newenvironment{proof}%
    {\rm \trivlist \item[\hskip \labelsep{\bf Proof. }]}%
    {\hspace*{\fill}$\Box$\endtrivlist}
    \newenvironment{varproof}%
    {\rm \trivlist \item[\hskip \labelsep{\bf Proof}]}%
    {\hspace*{\fill}$\Box$\endtrivlist}

\begin{document}
 \begin{center} \Large\bf
        A Christoffel Darboux formula for multiple orthogonal
        polynomials
    \end{center}

    \

    \begin{center}  \large
        E. Daems \\
            \normalsize \em
            Department of Mathematics, Katholieke Universiteit Leuven, \\
            Celestijnenlaan 200 B, 3001 Leuven, Belgium \\
            \rm evi.daems@wis.kuleuven.ac.be \\[3ex]
            \rm and \\[3ex]
            \large
        A.B.J. Kuijlaars\footnote{Supported by FWO Research Projects G.0176.02 and G.0455.04}  \\
            \normalsize \em
            Department of Mathematics, Katholieke Universiteit Leuven, \\
            Celestijnenlaan 200 B,  3001 Leuven, Belgium \\
            \rm arno@wis.kuleuven.ac.be
    \end{center}\ \\[1ex]

\begin{abstract}
Bleher and Kuijlaars recently showed that the eigenvalue
correlations from matrix ensembles with external source can
be expressed by means of a kernel built out of special multiple
orthogonal polynomials. We derive a Christoffel-Darboux formula for
this kernel for general multiple orthogonal polynomials. In
addition, we show that the formula can be written in terms of the
solution of the Riemann-Hilbert problem for multiple orthogonal
polynomials, which will be useful for asymptotic analysis.
\end{abstract}

\section{Introduction}
Multiple orthogonal polynomials are polynomials that satisfy orthogonal
conditions with respect to a number of weights, or more general with
respect to a number of measures. Such polynomials were first introduced
by Hermite in his proof of the transcendence of $e$, and were subsequently
used in number theory and approximation theory, see \cite{Apt}, \cite{ABVA},
and the references cited therein. The motivation for the present work comes
from a connection with random matrix theory. In the random matrix
model considered in \cite{BK1} the eigenvalue correlations are expressed
in terms of a kernel built out of multiple orthogonal polynomials with
respect to two weights
\begin{align} \label{twoweights}
    w_j(x) = e^{-V(x) +a_jx}, \qquad j=1,2, \qquad a_1 \neq a_2.
\end{align}
A Christoffel-Darboux formula was given in \cite{BK1} which leads to a
description of the kernel in terms of the Riemann-Hilbert
problem for multiple orthogonal polynomials \cite{VAGK}.
It is the aim of this paper to extend the Christoffel-Darboux formula
to multiple orthogonal polynomials with respect to an arbitrary number
of weights. We also allow more general weights than those in
(\ref{twoweights}).

Let $m \geq 2$ be an integer, and let
$w_1, w_2, \ldots, w_m$ be non-negative functions
on $\mathbb R$ such that all moments $\int_{-\infty}^{\infty} x^k w_j(x) dx$ exist.
Let $\vec{n} = (n_1, n_2, \ldots, n_m)$ be a vector of non-negative integers.
The (monic) multiple
orthogonal polynomial $P_{\vec{n}}$ of type II is a monic polynomial
of degree $|\vec{n}|$ satisfying
\begin{equation} \label{orco}
    \int P_{\vec{n}}(x) x^k w_j(x) dx = 0
    \qquad \mbox{for } k = 0, \ldots, n_j-1, \quad
        j = 1, \ldots, m.
\end{equation}
Here we define, as usual, $ |\vec{n}| = n_1 + n_2 + \cdots + n_m$.

We assume that the system is perfect, i.e., that for every
$\vec{n} \in (\mathbb N \cup \{0\})^m$, the polynomial
$P_{\vec{n}}$ exists and is unique, see \cite{Mah}. This is for
example the case when the weights form an Angelesco system or an
AT system, see e.g.\ \cite{VAC}.

The multiple orthogonal polynomials of type I are polynomials
$A_{\vec{n}}^{(k)}$ for $k=1,\ldots, m$, where $A_{\vec{n}}^{(k)}$
has degree $\leq n_k-1$, such that the function
\begin{align} \label{defQ}
Q_{\vec{n}}(x) = \sum_{k=1}^m A_{\vec{n}}^{(k)}(x) w_k(x)
\end{align}
 satisfies
\begin{align}\label{orthQ}
 \int x^j Q_{\vec{n}}(x) dx = \left\{
    \begin{array}{ll}
        0 & \mbox{ for } j = 0, \ldots, |\vec{n}| - 2, \\[10pt]
        1 & \mbox{ for } j= |\vec{n}| -1.
    \end{array} \right.
\end{align}
The polynomials $A_{\vec{n}}^{(k)}$ exist, are unique, and
they have full degree
\[ \deg A_{\vec{n}}^{(k)} = n_k-1, \]
since the system is perfect.

The usual monic orthogonal polynomials $P_n$ on the real
line with weight function $w(x)$ satisfy a three term recurrence
relation and this gives rise to the basic Christoffel-Darboux
formula
\begin{align}\label{christoffel}
\sum_{j=0}^{n-1}\frac{1}{h_j}P_j(x)P_j(y) =
    \frac{1}{h_{n-1}}\frac{P_{n}(x)P_{n-1}(y)-P_{n-1}(x)P_{n}(y)}{x-y},
\end{align}
where
\[ h_j = \int P_j(x) x^j w(x)dx. \]

In order to generalize the formula (\ref{christoffel}) to
multiple orthogonal polynomials, we consider a sequence
of multi-indices $\vec{n}_0, \vec{n}_1, \ldots, \vec{n}_n$ such that
for each $j = 0, 1, \ldots, n$,
\begin{align} \label{njs}
 |\vec{n}_j| = j, \qquad
\vec{n}_{j+1} \geq \vec{n}_j,
\end{align}
where the inequality is taken componentwise.
This means that we can go from $\vec{n}_j$ to $\vec{n}_{j+1}$
by increasing one of the components of $\vec{n}_j$ by $1$. We view
$\vec{n}_0, \vec{n}_1, \ldots, \vec{n}_n$ as a path from
$\vec{n}_0 = \vec{0}$ (the all-zero vector) to an arbitrary
multi-index $\vec{n} = \vec{n}_n$.
This path will be fixed
and all notions are related to this fixed path. Given such a
path, we define the polynomials $P_j$ and functions $Q_j$
(with single index) as
\begin{align}\label{naamPQ}
P_j = P_{\vec{n}_j}, \qquad Q_j = Q_{\vec{n}_{j+1}}.
\end{align}
Our aim is to find a simplified expression for the sum
\begin{align}\label{kern}
 K_n(x,y) = \sum_{j=0}^{n-1} P_j(x) Q_j(y).
\end{align}
To do this, we introduce the following notation. We define
for every multi-index $\vec{n}$ and every $k=1, \ldots, m$,
\begin{align}\label{defh}
h_{\vec{n}}^{(k)} = \int P_{\vec{n}}(x) x^{n_k} w_k(x) dx.
\end{align}
The numbers $h_{\vec{n}}^{(k)}$ are non-zero, since the
system is perfect.
We also use the standard basis vectors
\begin{align} \label{defen}
 \vec{e}_k = (0, \ldots, 0, 1, 0, \ldots, 0),
    \quad \mbox{where $1$ is in the $k$th position.}
\end{align}
Our result is the following.

\begin{theorem}\label{theorem}
Let $n \in \mathbb N$ and let $\vec{n}_0, \vec{n}_1, \ldots, \vec{n}_n$
be multi-indices such that {\rm (\ref{njs})} holds.
Let $P_j$ and $Q_j$ be as in {\rm (\ref{naamPQ})}.
Then we have if $\vec{n} = \vec{n}_n$,
\begin{align}\label{hoofd}
    (x-y)\sum_{j=0}^{n-1} P_j(x) Q_j(y) & =
    P_{\vec{n}}(x)Q_{\vec{n}}(y)-\sum_{k=1}^{m}
    \frac{h_{\vec{n}}^{(k)}}{h_{\vec{n} -\vec{e}_{k}}^{(k)}}
    P_{\vec{n}-\vec{e}_k}(x) Q_{\vec{n}+\vec{e}_k}(y).
\end{align}
\end{theorem}

It is easy to see that (\ref{hoofd}) reduces to the classical
Christoffel-Darboux formula (\ref{christoffel}) in case $m=1$. For
$m=2$ the formula was proven in \cite{BK1}.

\begin{remark}
It follows from (\ref{hoofd}) that the kernel (\ref{kern}) only
depends on the endpoint $\vec{n}$ of the chosen path from
$\vec{0}$ to $\vec{n}$ and not on the particular path itself,
since clearly the right-hand side of (\ref{hoofd}) only depends on
$\vec{n}$.

This fact can be deduced from the fact that
for any multi-index $\vec{k}$ and for $i \neq j$, we have
\begin{align} \label{relabel}
P_{\vec{k}}(x)Q_{\vec{k}+\vec{e_i}}(y) +
P_{\vec{k}+\vec{e}_i}(x)Q_{\vec{k}+\vec{e}_i+\vec{e}_j}(y) & =
P_{\vec{k}}(x)Q_{\vec{k}+\vec{e}_j}(y) +
P_{\vec{k}+\vec{e}_j}(x)Q_{\vec{k}+\vec{e}_i+\vec{e}_j}(y).
\end{align}
The relation (\ref{relabel}) follows easily from Lemma
\ref{theoremPQ}  below.
\end{remark}

\begin{remark}
In \cite{SVI}, Sorokin and Van Iseghem proved a Christoffel-Darboux
formula for vector polynomials that have matrix orthogonality
properties. As a special case this includes
the multiple orthogonal polynomials of type I and type II,
when one of the vector polynomials has only one component.
In this special case, their Christoffel-Darboux formula
comes down to the formula
\begin{align} \label{svi0}
    (x-y) \sum_{j=0}^{n-1} P_j(x) Q_j(y)
    = P_{n}(x) Q_{n-1}(y) - \sum_{k=n}^{n+m-1} \sum_{j=0}^{n-1}
        c_{j,k} P_j(x) Q_k(y)
\end{align}
where the constants $c_{j,k}$ are such that
\[ x P_k(x) = \sum_{j=0}^{k+1} c_{j,k} P_j(x), \]
see also (\ref{svi}) below. In the
setting of \cite{SVI} it holds that $c_{j,k} = 0$ if $k \geq j + m +1$,
so that the right-hand side
of (\ref{svi0}) has $1 + \frac{1}{2} m (m+1)$ terms. Note that in our formula
(\ref{hoofd}) the right-hand side has only $1+m$ terms.

Another Christoffel-Darboux formula for multiple orthogonal polynomials
similar to the one in \cite{SVI}  has been given recently in \cite{CVA}.
\end{remark}

\begin{remark}
As mentioned before, the formula (\ref{hoofd}) is useful in the
theory of random matrices. Br\'ezin and Hikami \cite{BH3} studied
a random matrix model with external source given by the
probability measure
\begin{align} \label{matrixmodel}
    \frac{1}{Z_n}e^{-Tr(V(M)-AM)}dM
\end{align}
defined on the space of $n \times n$ Hermitian matrices $M$. Here
we have that $V(x)=\frac{1}{2}x^2$,  $A$ is a fixed Hermitian
matrix (the external source), and $Z_n$ is a normalizing constant.
For this case, we can write $M = H+A$, where $H$ is a random
matrix from the Gaussian unitary ensemble and $A$ is
deterministic. Zinn-Justin \cite{ZJ} considered the case of an
arbitrary polynomial $V$.

The $k$-point correlation function $R_k(\lambda_1,\ldots,\lambda_k)$
of the (random) eigenvalues of a matrix from the ensemble
(\ref{matrixmodel}) can be expressed as a $k\times k$ determinant
involving a kernel $K_n(x,y)$
\begin{align}
R_k(\lambda_1,\ldots,\lambda_k) & =
\det(K_n(\lambda_i, \lambda_j))_{1\leq i,j \leq k},
\end{align}
see \cite{ZJ}. Suppose that the external source $A$ has $m$ distinct
eigenvalues $\alpha_1, \ldots, \alpha_m$ with respective
multiplicities $n_1, \ldots, n_m$.  Let $\vec{n} = (n_1, \ldots, n_m)$.
Then it was shown in \cite{BK1} that the
kernel $K_n$ has the form (\ref{kern}) built out of the multiple
orthogonal polynomials associated with the weights
\[ w_j(x) = e^{-(V(x) - \alpha_j x)}, \qquad j = 1, \ldots, m. \]
The Christoffel-Darboux formula (\ref{hoofd}) gives a compact expression
for the kernel.
\end{remark}

There is another expression for the kernel (\ref{kern}) in terms of
the solution of a Riemann-Hilbert problem. This will be especially useful
for the asymptotic analysis of the matrix model (\ref{matrixmodel}).
We will discuss this in the next section. The proof of
Theorem \ref{theorem} is presented in Section 3.

\section{Link with the Riemann-Hilbert problem}

Van Assche, Geronimo, and Kuijlaars \cite{VAGK} found a Riemann-Hilbert
problem that characterizes the multiple orthogonal polynomials. This is
an extension of the Riemann-Hilbert problem for orthogonal
polynomials due to Fokas, Its, and Kitaev \cite{FIK}.
We seek $Y:\mathbb{C} \setminus \mathbb{R} \to \mathbb{C}^{(m+1) \times (m+1)}$
such that
\begin{enumerate}
\item $Y$ is analytic on $\mathbb{C} \setminus \mathbb{R}$, \item
for $ x \in \mathbb{R}$, we have $Y_{+}(x)=Y_{-}(x)S(x)$, where
\begin{align}\label{RH1}
S(x) = \left[
\begin{array}{ccccc}
1 & w_1(x) & w_2(x) & \cdots & w_m(x) \\
0 & 1 & 0 & \cdots & 0 \\
0 & 0 & 1 & \cdots & 0 \\
\vdots & \vdots & \vdots & \ddots & \vdots \\
0 & 0 & 0 & \cdots & 1
\end{array} \right],
\end{align}
 \item as $z \to \infty$, we have
that
\begin{align}\label{RH11}
Y(z)=\left(I + O\left(\frac{1}{z}\right) \right) \left[
\begin{array}{ccccc}
z^n & 0 & 0 & \cdots & 0 \\
0 & z^{-n_1} & 0 &\cdots & 0 \\
0 & 0 & z^{-n_2} & \cdots & 0 \\
\vdots & \vdots & \vdots & \ddots & \vdots \\
0 & 0 & 0 & \cdots & z^{-n_m}
\end{array}
\right],
\end{align}
where $I$ denotes the $(m+1) \times (m+1)$ identity matrix.
\end{enumerate}
This Riemann-Hilbert problem has a unique solution given by:
\begin{align} \label{defY}
Y(z)=\left[
\begin{array}{cc}
P_{\vec{n}}(z) & \vec{R}_{\vec{n}}(z)\\
c_1 P_{\vec{n}-\vec{e_1}}(z) & c_1 \vec{R}_{\vec{n}-\vec{e}_1}(z)\\
c_2 P_{\vec{n}-\vec{e}_2}(z) & c_2 \vec{R}_{\vec{n}-\vec{e}_2}(z) \\
\vdots & \vdots \\
c_m P_{\vec{n}-\vec{e}_m}(z)& c_m \vec{R}_{\vec{n}-\vec{e}_m}(z)
\end{array}
\right],
\end{align}
where $P_{\vec{n}}(z)$ is the multiple orthogonal polynomial of
type II with respect to the weights $w_1, \ldots, w_m$ and
$\vec{R}_{\vec{n}}=(R_{\vec{n},1},R_{\vec{n},2},\ldots,R_{\vec{n},m})$
is the vector containing the Cauchy transforms
\[R_{\vec{n},j}(z)=\frac{1}{2 \pi i}\int \frac{P_{\vec{n}}(x)w_j(x)}{x-z}dx,\]
and
\begin{align} \label{defcj}
    c_j=-\frac{2 \pi i}{h^{(j)}_{\vec{n}-\vec{e}_j}}, \qquad j =1, \ldots, m.
\end{align}

Van Assche, Geronimo, and Kuijlaars \cite{VAGK} also gave a Riemann-Hilbert
problem that characterizes the multiple orthogonal polynomials of
type I. Here we seek $X: \mathbb{C} \setminus \mathbb{R} \to
\mathbb{C}^{(m+1)\times (m+1)}$ such that
\begin{enumerate}
\item $X$ is analytic on $\mathbb{C} \setminus \mathbb{R}$, \item
for $x \in \mathbb{R}$, we have $X_{+}(x)=X_{-}(x)U(x)$, where
\begin{align}\label{RH2}
U(x)=\left[
\begin{array}{ccccc}
1 & 0 & 0 & \cdots & 0\\
-w_1(x) & 1 & 0 & \cdots & 0 \\
-w_2(x) & 0 & 1 & \cdots & 0 \\
\vdots & \vdots & \vdots & \ddots & \vdots \\
-w_m(x) & 0 & 0 & \cdots & 1
\end{array}\right],
\end{align}
\item as $z \rightarrow \infty$, we have
\begin{align}\label{RH22}
X(z)=\left(I + O\left(\frac{1}{z}\right)\right) \left[
\begin{array}{ccccc}
z^{-n} & 0 & 0 & \cdots & 0 \\
0 & z^{n_1} & 0 & \cdots & 0 \\
0 & 0 & z^{n_2} & \cdots & 0 \\
\vdots & \vdots & \vdots & \ddots & \vdots \\
0 & 0 & 0 & \cdots & z^{n_m}
\end{array}
\right].
\end{align}
\end{enumerate}
This Riemann-Hilbert problem also has a unique solution and it is
given by
\begin{align} \label{defX}
    X(z)&=\left[
\begin{array}{cc}
\int Q_{\vec{n}}(x) \frac{dx}{z-x} & 2 \pi i \vec{A}_{\vec{n}}(z) \\
k_1 \frac{1}{2\pi i} \int Q_{\vec{n}+\vec{e}_1}(x) \frac{dx}{z-x} & k_1 \vec{A}_{\vec{n}+\vec{e}_1}(z) \\
k_2 \frac{1}{2\pi i} \int Q_{\vec{n}+\vec{e}_2}(x) \frac{dx}{z-x} & k_2 \vec{A}_{\vec{n}+\vec{e}_2}(z) \\
\vdots & \vdots \\
k_m \frac{1}{2\pi i} \int Q_{\vec{n}+\vec{e}_m}(x) \frac{dx}{z-x} & k_m \vec{A}_{\vec{n}+\vec{e}_m}(z)
\end{array}
\right],
\end{align}
where
$\vec{A}_{\vec{n}} = (A^{(1)}_{\vec{n}}, A^{(2)}_{\vec{n}}, \ldots, A^{(m)}_{\vec{n}})$
is the vector of multiple orthogonal polynomials of type I with respect
to $w_1, \ldots, w_m$, $Q_{\vec{n}}(z) = \sum_{k=1}^{m} A^{(k)}_{\vec{n}}(x)w_k(x)$
and
\begin{align} \label{defkj}
 k_j = h^{(j)}_{\vec{n}}, \qquad j =1, \ldots, m.
\end{align}

It is now possible to write the kernel $K_n(x,y)$ in terms of the
the solutions of the two Riemann-Hilbert problems, see also \cite{BK1}. First,
we observe that $X=Y^{-t}$. If we look at  the $j+1,1$-entry of the product
$Y^{-1}(y)Y(x)=X^t(y)Y(x)$, where $j=1,\ldots,m$, then we find by
(\ref{defY}) and (\ref{defX})
\begin{align} \nonumber
[Y^{-1}(y)Y(x)]_{j+1,1} &  =   \left[ \begin{array}{cccc} 2 \pi i
A^{(j)}_{\vec{n}}(y) & k_1 A^{(j)}_{\vec{n}+\vec{e}_1}(y) & \cdots
& k_m A^{(j)}_{\vec{n}+\vec{e}_m}(y) \end{array} \right] \left[
\begin{array}{c}
P_{\vec{n}}(x) \\
c_1 P_{\vec{n}-\vec{e_1}}(x) \\
c_2 P_{\vec{n}-\vec{e_2}}(x) \\
\vdots \\
c_m P_{\vec{n}-\vec{e_m}}(x)
\end{array}
\right] \\ \label{XtY} & =  2 \pi i
\left(P_{\vec{n}}(x)A^{(j)}_{\vec{n}}(y) - \sum_{k=1}^m
\frac{h^{(k)}_{\vec{n}}}{h^{(k)}_{\vec{n}-\vec{e}_k}}
P_{\vec{n}-\vec{e}_k}(x)A^{(j)}_{\vec{n}+\vec{e}_k}(y)\right).
\end{align}
where in the last step we used the expressions (\ref{defcj}) and
(\ref{defkj}) for the constants $c_j$ and $k_j$.
Multiplying (\ref{XtY}) by $w_j(y)$, dividing by $2 \pi i$, and
summing over $j=1, \ldots, m$, we obtain the right-hand side of
(\ref{hoofd}). Therefore we see that
\begin{align} \nonumber
(x-y) K_n(x,y) & =
    \frac{1}{2\pi i} \sum_{j=1}^{m} w_j(y)[Y^{-1}(y)Y(x)]_{j+1,1}\\
    \label{hoofd2}
    & =\frac{1}{2\pi i} \left[
\begin{array}{cccc}
0 & w_1(y)  & \cdots & w_m(y)
\end{array}\right]
Y^{-1}(y)Y(x) \left[
\begin{array}{c}
1 \\ 0 \\ \vdots \\ 0
\end{array}
\right].
\end{align}
It is clear that the right-hand side of (\ref{hoofd2}) is $0$ for $x=y$,
which is not obvious at all for the right-hand side of (\ref{hoofd}).

In \cite{BK2} the Riemann-Hilbert problem (\ref{RH1})--(\ref{RH11})
is analyzed in the limit $n \to \infty$ for the special case of $m=2$,
$n_1=n_2$, and weights
\[ w_1(x) = e^{- n(\frac{1}{2} x^2 - ax)}, \qquad w_2(x) = e^{-n(\frac{1}{2}x^2 + ax)}. \]
The corresponding multiple orthogonal polynomials are known as multiple
Hermite polynomials \cite{ABVA,VAC}. The asymptotic analysis of (\ref{RH1})--(\ref{RH11})
is done by the Deift/Zhou steepest descent method for Riemann-Hilbert problems,
see \cite{Dei} and references cited therein.

\section{Proof of Theorem \ref{theorem}}
For the proof we are going to extend the path $\vec{n}_0, \vec{n}_1,
\ldots, \vec{n}_n$ by defining
\begin{align}\label{assume}
 \vec{n}_{n+k} - \vec{n}_{n+k-1} = \vec{e}_k,
    \qquad k = 1, 2, \ldots, m.
\end{align}
We will also extend the definition (\ref{naamPQ}) by
putting $P_j = P_{\vec{n}_j}$ and $Q_{j-1} = Q_{\vec{n}_j}$ for
$j = n+1, \ldots, n+m$.

\subsection{Biorthogonality and recurrence relations}
The multiple orthogonal polynomials satisfy a
biorthogonality relation.
\begin{lemma}\label{biorth}
    We have
\[ \int P_k(x) Q_j(x) dx = \delta_{j,k}, \]
    where $\delta_{j,k}$ is the Kronecker delta.
\end{lemma}
\begin{proof}
This is immediate from the definitions (\ref{naamPQ}), the
orthogonality conditions (\ref{orthQ}) of the function $Q_j$ and
(\ref{orco}) of the polynomial $P_k$ and the fact that $P_k$ is a
monic polynomial.
\end{proof}

Because $xP_k(x)$ is a polynomial of degree $k+1$, we can expand
$xP_k(x)$ as
\begin{align}\label{expP}
 x P_k(x) & = \sum_{j=0}^{k+1} c_{j,k} P_j(x).
\end{align}
The coefficients can be calculated by Lemma \ref{biorth} by
multiplying both sides of (\ref{expP}) with $Q_j(x)$ and
integrating over the real line. That gives us
\begin{align} \label{defcjk}
    c_{j,k}=\int x P_k(x) Q_j(x) dx.
\end{align}
The coefficients $c_{j,k}$ are $0$ if $j \geq k+2$.

Because of (\ref{assume}) we can write $yQ_j(y)$ with $j \leq n-1$ as a linear
combination of $Q_0, \ldots, Q_{n+m-1}$ and we have by Lemma \ref{biorth}
\begin{align}\label{expQ}
yQ_j(y) = \sum_{k=0}^{n+m-1} c_{j,k} Q_k(y) \qquad
    \mbox{for } j=0,\ldots,n-1.
\end{align}

Using the expansions (\ref{expP}) and (\ref{expQ}) for $xP_k(x)$
and $y Q_j(y)$ we can write
\begin{align*}
(x-y)\sum_{k=0}^{n-1}P_k(x)Q_k(y) & =
    \sum_{k=0}^{n-1}xP_k(x)Q_k(y)-\sum_{k=0}^{n-1}P_k(x)yQ_k(y) \nonumber \\
    & =  \sum_{k=0}^{n-1}\sum_{j=0}^{k+1}
    c_{j,k}P_j(x)Q_k(y)-\sum_{k=0}^{n-1}
    \sum_{j=0}^{n+m-1}c_{k,j}P_k(x)Q_j(y).
\end{align*}
A lot of terms cancel. Since $c_{j,k} = 0$ for $j \geq k+2$,
and $c_{n,n-1} = 1$, what remains is
\begin{align}\label{svi}
(x-y) \sum_{k=0}^{n-1} P_k(x) Q_k(y) & =
    P_{\vec{n}}(x)Q_{\vec{n}}(y) -
    \sum_{k=n}^{n+m-1} \sum_{j=0}^{n-1} c_{j,k} P_j(x) Q_k(y).
\end{align}
We also used the fact that $P_n = P_{\vec{n}}$ and $Q_{n-1} = Q_{\vec{n}}$.
Note that (\ref{svi}) corresponds to the Christoffel-Darboux
formula of \cite{SVI}, as mentioned in the introduction.

In the rest of the proof we are going to show that
\begin{align} \label{tbw}
    \sum_{k=n}^{n+m-1} \sum_{j=0}^{n-1} c_{j,k} P_j(x) Q_k(y)
    = \sum_{k=1}^m \frac{h_{\vec{n}}^{(k)}}{h_{\vec{n}-\vec{e}_k}^{(k)}}
        P_{\vec{n}-\vec{e}_k}(x) Q_{\vec{n}+\vec{e}_k}(y)
\end{align}
so that (\ref{svi}) then leads to our desired formula (\ref{hoofd}).

\subsection{The vector space generated by the polynomials
$P_{\vec{n}-\vec{e}_1},\ldots,P_{\vec{n}-\vec{e}_m}$}

For fixed $y$, the right-hand side of (\ref{tbw}) belongs to the vector
space spanned by the polynomials
of $P_{\vec{n}-\vec{e}_1}, \ldots, P_{\vec{n}-\vec{e}_m}$.
In this part of the proof, we characterize this vector space
and show that the left-hand side of (\ref{tbw}) also belongs
to this vector space $V$.

\begin{lemma}\label{vector}
The polynomials
$P_{\vec{n}-\vec{e}_1},\ldots,P_{\vec{n}-\vec{e}_m}$ are a basis
of the vector space $V$ of all polynomials $\pi$ of degree $\leq n-1$ satisfying
\begin{align}\label{dependent}
     \int \pi(x) x^{i} w_j(x) dx = 0,
    \qquad i=0,\ldots, n_j-2, \quad j=1,\ldots,m.
\end{align}
\end{lemma}
\begin{proof}
By the orthogonality properties (\ref{orco}) of the polynomials
$P_{\vec{n}-\vec{e}_i}$ for $i =1, \ldots, m$, it is obvious that
they belong to $V$. We are first going to show that the
polynomials $P_{\vec{n}-\vec{e}_i}$ are linearly independent.
Suppose that
\begin{align}\label{linonafh}
a_1P_{\vec{n}-\vec{e}_1} + a_2P_{\vec{n}-\vec{e}_2} + \cdots + a_m
P_{\vec{n}-\vec{e}_m} & = 0
\end{align}
for some coefficients $a_j$.  Multiplying (\ref{linonafh}) with
$w_j(x)x^{n_j-1}$, and integrating over the real line, we obtain
$ a_j h^{(j)}_{\vec{n}-\vec{e}_j}=0$. Since $h^{(j)}_{\vec{n}-\vec{e}_j} \neq 0$,
we get $a_j =0$ for $j = 1, \ldots, m$, which shows that the polynomials
are linearly independent.

Suppose next that $\pi$ belongs to $V$. Put
\[ b_j = \frac{1}{h^{(j)}_{\vec{n}-\vec{e}_j}} \int \pi(x) x^{n_j-1} w_j(x) dx \]
and define the polynomial $\pi_1$ by
\begin{align} \label{defpi1}
\pi_1 & = b_1P_{\vec{n}-\vec{e}_1} + b_2P_{\vec{n}-\vec{e}_2}
+\cdots + b_m P_{\vec{n}-\vec{e}_m}.
\end{align}
Then $\pi_1 - \pi$ belongs to $V$ and
\begin{align}
\int \left(\pi_1(x) - \pi(x) \right) x^{n_j-1} w_j(x)dx & = 0,
    \qquad j = 1, \ldots, m.
\end{align}
This means that $\pi_1-\pi$ satisfies the conditions
\begin{align} \label{piispi1}
     \int \left(\pi_1(x) - \pi(x) \right) x^{i} w_j(x) dx = 0,
    \qquad i=0, \ldots, n_j-1, \quad j=1,\ldots, m.
\end{align}
Because $\pi_1-\pi$ is a polynomial of degree $\leq n-1$ and the
system is perfect, it follows from (\ref{piispi1}) that $\pi_1 - \pi = 0$. Therefore
$\pi=\pi_1$, and $\pi$ can be written as a linear combination of
the polynomials $P_{\vec{n}-\vec{e}_1},\ldots,P_{\vec{n}-\vec{e}_m}$.

The lemma follows.
\end{proof}

\begin{lemma} \label{pikinV}
For every $k = n, \ldots, n+m-1$, we
have that the polynomial
\begin{align} \label{defpik}
    \pi_k(x) = \sum_{j=0}^{n-1} c_{j,k} P_j(x)
\end{align}
belongs to the vector space $V$.
\end{lemma}

\begin{proof}
Clearly $\pi_k$ is a polynomial of degree $n-1$.
Using (\ref{expP}) we see that
\begin{align} \label{representatie}
 \pi_k(x) = x P_k(x) - \sum_{j=n}^{k+1} c_{j,k} P_j(x).
\end{align}
The representation (\ref{representatie}) of $\pi_k$
and the orthogonality conditions (\ref{orco})
show that $\pi_k$ satisfies the relations (\ref{dependent}),
so that $\pi_k$ belongs to $V$ by Lemma \ref{vector}.
\end{proof}

Because of Lemma \ref{pikinV}, the left-hand side of (\ref{tbw})
belongs to $V$ for every $y$, and so by Lemma \ref{vector},
we can write
\begin{align} \label{phikes}
    \sum_{k=n}^{n+m-1} \sum_{j=0}^{n-1} c_{j,k} P_j(x) Q_k(y)
       % = \sum_{k=n}^{n+m-1} \sum_{j=1}^m b_{j,k} P_{\vec{n}-\vec{e}_j}(x)
       %     Q_k(y)
        = \sum_{j=1}^m \phi_j(y) P_{\vec{n}-\vec{e}_j}(x)
\end{align}
for certain functions $\phi_j(y)$. The next lemma gives
an expression for $\phi_j$.
We use the notation $\vec{s}_0 = \vec{0}$ (all-zero vector) and
\[ \vec{s}_j = \sum_{k=1}^j \vec{e}_k, \qquad j = 1, \ldots, m. \]

\begin{lemma}
We have for $j=1, \ldots, m$,
\begin{align} \label{phijexpr}
    h_{\vec{n}-\vec{e}_j}^{(j)} \phi_j(y) =
    \sum_{i=1}^{j} h_{\vec{n} + \vec{s}_{i-1}}^{(j)} Q_{\vec{n}+\vec{s}_i}(y).
\end{align}
\end{lemma}
\begin{proof}
Rewriting the left-hand side of (\ref{phikes}) using
(\ref{defpik}) and (\ref{representatie}) we obtain
\begin{align} \label{phikes2}
\sum_{j=1}^m \phi_j(y) P_{\vec{n}-\vec{e}_j}(x)
    = \sum_{k=n}^{n+m-1} x P_k(x) Q_k(y)
        - \sum_{k=n}^{n+m-1} \sum_{j=n}^{k+1} c_{j,k} P_j(x) Q_k(y).
\end{align}

Now multiply (\ref{phikes2}) with $x^{n_j-1} w_j(x)$ and
integrate with respect to $x$. Then the left-hand
side gives
\begin{align}\label{leftside}
h_{\vec{n}-\vec{e}_j}^{(j)} \phi_j(y).
\end{align}
The second sum in the right-hand side of (\ref{phikes2}) gives no
contribution to the integral because of orthogonality,
and the first sum gives
\begin{align} \label{rightside}
 \sum_{k=n}^{n+m-1} \left( \int P_k(x) x^{n_j}
w_j(x) dx \right) Q_k(y) =
    \sum_{i=1}^m
    \left( \int P_{n+i-1}(x) x^{n_j} w_j(x) dx \right)
        Q_{n+i-1}(y).
\end{align}
Because of the choice (\ref{assume}) and the definition (\ref{naamPQ})
we have
\[ P_{n+i-1} = P_{\vec{n}+ \vec{s}_{i-1}},
    \quad Q_{n+i-1} = Q_{\vec{n} + \vec{s}_i}
    \qquad  i = 1, \ldots, m. \]
Then we see that the integral in the right-hand side of
(\ref{rightside}) is zero if $i \geq j+1$ and otherwise it is
equal to $h_{\vec{n} + \vec{s}_{i-1}}^{(j)}$. Then
(\ref{phijexpr}) follows.
\end{proof}

\subsection{Completion of the proof of Theorem \ref{theorem}}

In view of (\ref{phikes}) and
(\ref{phijexpr}) it remains to prove that
\begin{align} \label{tbw2}
    h_{\vec{n}}^{(j)} Q_{\vec{n}+\vec{e}_j} =
    \sum_{i=1}^j h_{\vec{n} + \vec{s}_{i-1}}^{(j)}
        Q_{\vec{n}+\vec{s}_i}(y)
\end{align}
for $j=1, \ldots, m$, and then (\ref{tbw}) follows.

To establish (\ref{tbw2}) we need some properties of
the numbers $h_{\vec{n}}^{(j)}$ and relations between
$Q$-functions with different multi-indices.
We already noted that $h_{\vec{n}}^{(j)} \neq 0$.
We express the leading coefficients of the polynomials
$A^{(j)}_{\vec{n}}$ in terms of these numbers.
\begin{lemma}\label{theoremlc}
The leading coefficient of
$A_{\vec{n}+\vec{e}_j}^{(j)}$ is equal to
$\frac{1}{h_{\vec{n}}^{(j)}}$.
\end{lemma}

\begin{proof}
Because of the orthogonality conditions (\ref{orco}) and
(\ref{orthQ}) we have that
\begin{eqnarray*}
1 & = & \int P_{\vec{n}}(x)Q_{\vec{n}+\vec{e}_j}(x) dx \\
& = & \int
P_{\vec{n}}(x)A_{\vec{n}+\vec{e}_j}^{(j)}(x) w_j(x)dx \\
& =& \left( \mbox{leading coefficient of } A_{\vec{n}+\vec{e}_j}^{(j)} \right)
    \int P_{\vec{n}}(x)x^{n_j}w_j(x)dx \\
& = & \left( \mbox{leading coefficient of } A_{\vec{n}+\vec{e}_j}^{(j)} \right) h_{\vec{n}}^{(j)}
\end{eqnarray*}
and the lemma follows.
\end{proof}

\begin{lemma}\label{theoremPQ}
Let $j \neq k$. Then we have for every multi-index $\vec{n}$ that
\begin{align} \label{contigP}
P_{\vec{n}}(x) & =
    \frac{h_{\vec{n}}^{(k)}}{h_{\vec{n}+\vec{e}_j}^{(k)}}
    \left(P_{\vec{n}+\vec{e}_j} - P_{\vec{n}+\vec{e}_k}\right)
    =
    - \frac{h_{\vec{n}}^{(j)}}{h_{\vec{n}+\vec{e}_k}^{(j)}}
    \left(P_{\vec{n}+\vec{e}_j} - P_{\vec{n}+\vec{e}_k} \right)
\end{align}
and
\begin{align} \label{contigQ}
Q_{\vec{n}}& =
\frac{h_{\vec{n}-\vec{e}_j - \vec{e}_k}^{(k)}}{h_{\vec{n}-\vec{e}_k}^{(k)}}
\left(Q_{\vec{n} - \vec{e}_j} -Q_{\vec{n}-\vec{e}_k}\right)
=
-\frac{h_{\vec{n}-\vec{e}_j - \vec{e}_k}^{(j)}}{h_{\vec{n}-\vec{e}_{j}}^{(j)}}
\left(Q_{\vec{n}-\vec{e}_j}-Q_{\vec{n}-\vec{e}_k}\right).
\end{align}
\end{lemma}
\begin{proof}
We know that $P_{\vec{n}}$ is a polynomial of degree $|\vec{n}|$
that satisfies the orthogonality conditions (\ref{orco}). It is
easy to see that $P_{\vec{n}+\vec{e}_j} - P_{\vec{n}+\vec{e}_k}$
is a polynomial of degree $|\vec{n}|$ that satisfies these same
conditions. Because the system is perfect, we then have that
\begin{align}\label{bewijs1}
\gamma P_{\vec{n}}(x) & = P_{\vec{n}+\vec{e}_j}(x)-P_{\vec{n} +
\vec{e_k}}(x),
\end{align}
for some $\gamma \in \mathbb{R}$.
Multiplying (\ref{bewijs1}) with $x^{n_k} w_k(x)$ and
integrating over the real line, we find that
\[ \gamma h_{\vec{n}}^{(k)} = h_{\vec{n} +\vec{e_j}}^{(k)} - 0
    = h_{\vec{n} + \vec{e_j}}^{(k)}.
    \]
This proves the first equality of (\ref{contigP}). The second equality
follows by interchanging $j$ and $k$.
\medskip

Next we show (\ref{contigQ}).
It is easy to see that $Q_{\vec{n}-\vec{e}_j}- Q_{\vec{n}-\vec{e}_k}$
satisfies the same orthogonality conditions (\ref{orthQ}) as $Q_{\vec{n}}$.
Since the degrees of the polynomials $A_{\vec{n}-\vec{e}_j}^{(i)} -
A_{\vec{n}-\vec{e}_k}^{(i)}$ do not exceed the degrees of
$A_{\vec{n}}^{(i)}$ for $i=1, \ldots, m$, it follows that
\begin{align}\label{bewijs2}
\gamma Q_{\vec{n}} & = Q_{\vec{n}-\vec{e}_j}-Q_{\vec{n}-\vec{e}_k},
\end{align}
for some $\gamma \in \mathbb{R}$. To compute $\gamma$, we are
going to compare the leading coefficients of the polynomials that
come with $w_k(x)$. Using Lemma \ref{theoremlc}, we find
that
\[ \gamma \frac{1}{h_{\vec{n}-\vec{e}_k}^{(k)}} =
    \frac{1}{h_{\vec{n}-\vec{e}_j-\vec{e}_k}^{(k)}} - 0  =
    \frac{1}{h_{\vec{n}-\vec{e}_j-\vec{e}_k}^{(k)}}. \]
This proves the first equality of (\ref{contigQ}). The second equality follows
by interchanging $j$ and $k$.
\end{proof}

Now we are ready to complete  the proof of Theorem \ref{theorem}.
\begin{varproof} \textbf{of Theorem \ref{theorem}}

In view of what was said before, it suffices to prove (\ref{tbw2}).
Fix $j=1, \ldots, m$.
We are going to prove by induction that for $k = 0, \ldots, j-1$,
\begin{align} \label{eqQ}
 h_{\vec{n}}^{(j)} Q_{\vec{n}+\vec{e}_j} =
 \sum_{i=1}^k h_{\vec{n}+\vec{s}_{i-1}}^{(j)}
 Q_{\vec{n}+\vec{s}_i}(y)
 + h_{\vec{n}+\vec{s}_k}^{(j)} Q_{\vec{n}+\vec{s}_k+\vec{e}_j}.
\end{align}
For $k=0$, the sum in the right-hand side of (\ref{eqQ}) is an empty sum,
and then the equality (\ref{eqQ}) is clear.

Suppose that (\ref{eqQ}) holds for some $k \leq j-2$.  Taking
(\ref{contigQ}) with $\vec{n} + \vec{s}_{k+1}+\vec{e}_j$
instead of $\vec{n}$ and $k+1$ instead of $k$, we get
\begin{align*}
 Q_{\vec{n}+\vec{s}_{k+1}+\vec{e}_j} & =
-\frac{h^{(j)}_{\vec{n}+\vec{s}_{k}}}
{h_{\vec{n} + \vec{s}_{k+1}}^{(j)}}
\left(Q_{\vec{n}+\vec{s}_{k+1}}- Q_{\vec{n}+\vec{s}_k+\vec{e}_j}\right).
\end{align*}
Thus
\begin{align} \label{tussenstap}
h^{(j)}_{\vec{n}+\vec{s}_{k}} Q_{\vec{n}+\vec{s}_{k}+\vec{e}_j} & =
h_{\vec{n} + \vec{s}_{k}}^{(j)} Q_{\vec{n}+\vec{s}_{k+1}}
+ h^{(j)}_{\vec{n}+\vec{s}_{k+1}} Q_{\vec{n}+\vec{s}_{k+1}+\vec{e}_j}
\end{align}
and using  the induction hypothesis (\ref{eqQ}) we obtain
(\ref{eqQ}) with $k$ replaced by $k+1$.

So (\ref{eqQ}) holds for every $k = 0, 1, \ldots, j-1$. Taking
$k=j-1$ in (\ref{eqQ}), we obtain (\ref{tbw2}) and this completes
the proof of Theorem \ref{theorem}.
\end{varproof}

\end{document}